\documentclass[12pt,reqno]{amsart}
\newtheorem{thm}{Theorem}[section]
\newtheorem{defi}{Definition}[section]

\newtheorem{cor}{Corollary}[section]
\newtheorem{pr}{Proposition}[section]
\theoremstyle{definition}
\newtheorem{rem}{Remark}[section]

\newcommand{\be}{\begin{equation}}
\newcommand{\ee}{\end{equation}}
\newcommand{\bea}{\begin{eqnarray}}
\newcommand{\eea}{\end{eqnarray}}
\newcommand{\beb}{\begin{eqnarray*}}
\newcommand{\eeb}{\end{eqnarray*}}
\usepackage{amssymb,amsfonts,amsthm,setspace,indentfirst,mathrsfs}
\usepackage{minibox}
\usepackage{enumitem}
\usepackage{cite}
\numberwithin{equation}{section}
\begin{document}
%%%%%%%%%%%%%%%%%%%%%%%%%%%%%%%%%%%%%%%%%%%%%%%%%%%%%%%%%%%%%%%%%%%%%%%%%%%%%%%%%%%%%%%%%%%%%%%%%%%%%%%%%%%%%
%
\title{On curvature related geometric properties of Hayward black hole spacetime}
\author[A. A. Shaikh., S. K. Hui., B. R. Datta and M. Sarkar]{Absos Ali Shaikh$^{* 1}$, Shyamal kumar Hui$^2$, Biswa Ranjan Datta$^3$ and Mousumi Sarkar$^4$}
\date{\today}
\address{\noindent\newline$^{1,2,3,4}$ Department of Mathematics,
	\newline University of Burdwan, 
	\newline Golapbag, Burdwan-713104,
	\newline West Bengal, India} 
\email{aask2003@yahoo.co.in$^1$, aashaikh@math.buruniv.ac.in$^1$}
\email{skhui@math.buruniv.ac.in$^2$}
\email{biswaranjandatta2019@gmail.com$^3$}
\email{mousrkr13@gmail.com$^4$}

\dedicatory{}
%%%%%%%%%%%%%%%%%%%%%%%%%%%%%%%%%%%%%%%%%%%%% Abstract %%%%%%%%%%%%%%%%%%%%%%%%%%%%%%%%%%%%%%%%%%%%%%%%%%%%%

\begin{abstract}
	
	This paper is devoted to the study of curvature properties  of Hayward black hole (briefly, HBH) spacetime, which is a solution of Einstein field equations (briefly, EFE) having non-vanishing cosmological constant. We have proved that the HBH spacetime is an Einstein manifold of level $2$, $2$-quasi Einstein, generalized quasi-Einstein and  Roter type manifold.  Also, it is shown that the nature of the HBH spacetime is pseudosymmetric and it obeys  several types of pseudosymmetries, such as,  pseudosymmetry due to concircular, conformal and conharmonic curvature (i.e., $F\cdot F=\mathcal{L}Q(g,F)$   for $F=W,C, K$ with a smooth scalar function $ \mathcal{L} $), and it also possesses the relation $R\cdot R-\mathcal{L} Q(g,C)=Q(S,R)$. It is engrossing to mention that the nature of energy momentum tensor of the HBH spacetime is pseudosymmetric.	On the basis of curvature related properties, we have made a comparison among Reissner-Nordstr\"om spacetime, interior black hole spacetime and HBH spacetime. Also,  it is shown that the HBH spacetime admits an almost $\eta$-Ricci soliton as well as an almost $\eta$-Ricci-Yamabe soliton. Finally, an elegant comparative study is delineated between the HBH spacetime and the point-like global monopole spacetime with respect to different kinds of symmetry, such as, motion, curvature collineation, curvature inheritance etc.

\end{abstract}

%%%%%%%%%%%%%%%%%%%%%%%%%%%%%%%%%%%%%%%%%%%%%%%%%%%%%%%%%%%%%%%%%%%%%%%%%%%%%%%%%%%%%%%%%%%%%%%%%%%%%%%%%%%%
%
\maketitle
%

%%%%%%%%%%%%%%%%%%%%%%%%%%%%%%%%%%%%%%%%%%%%%%%%%%%%%%%%%%%%%%%%%%%%%%%%%%%%%%%%%%%%%%%%%%%%%%%%%%%%%%%%%%%%%%
%																					Introduction
%%%%%%%%%%%%%%%%%%%%%%%%%%%%%%%%%%%%%%%%%%%%%%%%%%%%%%%%%%%%%%%%%%%%%%%%%%%%%%%%%%%%%%%%%%%%%%%%%%%%%%%%%%%%%%
\section{\bf Introduction}\label{intro}
%%%%%%%%%%%%%%%%%%%%%%%%%%%%%%%%%%%%%%%

Let us consider a semi-Riemannian manifold $M$ of dimension $n\geq 3$ such that $\nabla$ is the Levi-Civita connection of the semi-Riemannian metric $g$ with signature $(t,n-t)$, $0\leq t\leq n$ and $R$, $S$, $\kappa$ are respectively the Riemann, Ricci, scalar curvature of $M$. A connected 4-dimensional manifold $M$ with Lorentzian signature (1,3) or (3,1) is a spacetime.\\

\indent The curvature carries an enormous significance to acquire the shape of a space. In fact, the geometry of a space can be described by curvature explicitly as the relation  $\nabla R=0$ defines the notion of locally symmetric manifolds  (see, Cartan \cite{Cart26}). Cartan \cite{Cart46} introduced the notion of semisymmetric manifolds defined as $R\cdot R=0$ (see also, \cite{Szab82,Szab84,Szab85}) and the concept of pseudosymmetric manifolds was introduced by Adam\'{o}w and Deszcz \cite{AD83}, which is known as Deszcz pseudosymmetric space. A large number of physicists and mathematicians  investigated the concept of locally symmetric manifolds and introduced several generalized  notions of symmetries, such as, recurrent manifolds by Ruse \cite{Ruse46, Ruse49a, Ruse49b}, (see also \cite{Walk50}), different kinds of generalized notion of recurrent manifolds by Shaikh et al. \cite{SP10,SR10,SRK15,SRK18,SAR13,SR11,SKA18,SRK17}, curvature $2$-forms of recurrent manifolds by Besse \cite{Bess87, LR89, MS12a,MS13a,MS14}, pseudosymmetric manifolds by Chaki \cite{Chak87, Chak88}, weakly symmetric manifolds by Tam\'{a}ssy and Binh \cite{TB89, TB93} etc.  Haesen and Verstraelen \cite{HV07,HV07A,HV09} exhibited the geometrical and physical significance of various pseudosymmetries. We mention that the   Deszcz pseudosymmetry  achieved a great importance during last four decades due to its applications in the study of general relativity and cosmology  as numerous spacetimes (see, \cite{ ADEHM14, DK99, Kowa06, SAA18, SAAC20,  SDKC19, SAACD_LTB_2022,EDS_sultana_2022}) have been found to be pseudosymmetric. It is noteworthy  to mention that pseudosymmetries in the sense of Deszcz  and  Chaki  are not equivalent (see, \cite{SDHJK15}). \\

In 1982, during the study of compact $3$-dimensional manifolds with positive Ricci curvature, Hamilton \cite{Hamilton1982} established a process of evolving a Riemanninan metric over time, called Ricci flow. The self-similar solutions of the Ricci flow are known as Ricci solitons, which are natural generalizations of Einstein metrics \cite{Bess87,Brink1925,SYH09,S09}. The notion of Ricci soliton has been generalized in different ways, e.g. almost Ricci soliton,   $\eta$-Ricci soliton, almost $\eta$-Ricci soliton etc.

If the Ricci curvature $S$ and the metric tensor $g$ of a Riemannian manifold $M$ realize 
$$\frac{1}{2}\pounds_\xi g+S-\mu g=0$$
for a constant $\mu$, then $M$ is said to be a Ricci soliton, where $\pounds_\xi$ is the Lie derivative in the direction of the soliton vector field $\xi$. It is expanding, steady or shrinking  according to the condition $\mu<0$, $\mu=0$ or $\mu>0$ respectively. It is called an almost Ricci soliton \cite{Pigola2011} if $\mu$ is a non-constant  smooth function. We mention that if the corresponding soliton vector field $\xi$ of a Ricci soliton is Killing, then the Ricci soliton turns into an Einstein manifold. Again, if  a non-zero $1$-form $\eta$ on $M$ satisfies the relation 
$$\frac{1}{2}\pounds_\xi g+S-\mu g+\lambda (\eta\otimes\eta)=0,$$
$\mu,\lambda$ being constants, then $M$ is called an $\eta$-Ricci soliton \cite{Cho2009}. The $\eta$-Ricci soliton is said to be an almost $\eta$-Ricci soliton \cite{Blaga2016} if $\mu, \lambda$ are allowed to be smooth functions. \\

On the other hand, simultaneously with the notion of Ricci flow, Hamilton \cite{Hamilton1988} introduced the notion of Yamabe flow. Recently, as a scalar combination of Ricci and Yamabe flow, G\"uler and Cr\'{a}\c{s}mare\v{a}nu \cite{Guler2019} established a new geometric flow, which is
called Ricci-Yamabe flow, and Ricci-Yamabe  (resp., Yamabe) solitions are the self-similar solutions of Ricci-Yamabe  (resp., Yamabe) flow. If in a Riemannian manifold $M$  the Ricci curvature $S$ and the metric tensor $g$ realize the relation 
$$\frac{1}{2}\pounds_\xi g+\alpha_1 S+\left(\mu-\frac{1}{2}\alpha_2\kappa\right) g=0,$$
with the constants $\alpha_1$, $\alpha_2$, $\mu$, scalar curvature $\kappa$ and the soliton vector field $\xi$, then $M$ is called a Ricci-Yamabe soliton \cite{Siddiqi2020}. We note that if $\alpha_1=0,\alpha_2=2$  (resp.,  $\alpha_1=1,\alpha_2=0$), then it turns into Yamabe soliton (resp., Ricci soliton). In addition, if $\alpha_1$, $\alpha_2$, $\mu$ are allowed to be non-constant smooth functions, then $M$ is known as an almost Ricci-Yamabe soliton \cite{Siddiqi2020}. Again, if there is a non-zero 1-form $\eta$ satisfying
$$\frac{1}{2}\pounds_\xi g+\alpha_1 S+\left(\mu-\frac{1}{2}\alpha_2\kappa\right) g+\lambda \eta\otimes\eta=0,$$
with the constants $\alpha_1$, $\alpha_2$, $\mu$, $\lambda$, then $M$ is called an   $\eta$-Ricci-Yamabe soliton \cite{Siddiqi2020}. If the constants $\alpha_1$, $\alpha_2$, $\mu$, $\lambda$ are allowed to be non-constant smooth functions, then $M$ is called an almost $\eta$-Ricci-Yamabe soliton \cite{Siddiqi2020}.   A plenty of research papers (see, \cite{AliAhsan2013, AliAhsan2015, Ahsan2018,SDAA_LCS_2021} and the references therein) on Ricci soliton, Yamabe soliton and their generalizations are appeared during last three decades, and  nowadays it is an abuzz topic of research in differential geometry.\\

To construct gravitational potential, one can impose the symmetry in EFE and hence the geometrical symmetries play a crucial role in the theory of general relativity. Along a vector field, certain geometric quantity is preserved  if the Lie derivative of the corresponding tensor vanishes with respect to that vector field, and the vanishing Lie derivative explains the geometrical symmetries. The notions of motion, curvature collineation, Ricci collineation etc. are  the examples of  such symmetries. Katzin {et al.} \cite{KLD1969,KLD1970} rigorously investigated the role of curvature collineation in general relativity. In 1992, Duggal \cite{Duggal1992}  introduced the notion of curvature inheritance generalizing the concept of curvature collineation for the (1,3)-type curvature tensor. During last three decades, a plenty of articles (see, \cite{Ahsan1978,Ahsan1977_231,Ahsan1977_1055,Ahsan1987,Ahsan1995,Ahsan1996,Ahasan2005,AhsanAli2014,AA2012, AH1980,AliAhsan2012,SASZ2022,ShaikhDatta2022}) 
appeared in the literature regarding the investigations of such kinds of symmetries. 
Recently, during the investigation of geometric properties of Robinson-Trautman spacetime, Shaikh and Datta \cite{ShaikhDatta2022} introduced the concept of generalized curvature inheritance, which is a generalization of curvature collineation as well as curvature inheritance for the (0,4)-type curvature tensor. We note that the notions of curvature inheritance for the (1,3)-type curvature tensor and for the (0,4)-type curvature tensor are not equivalent \cite{ShaikhDatta2022}. In this paper, we have checked that the HBH spacetime does not admit any of the curvature related symmetries. Finally,  a worthy comparison between the HBH spacetime and the point-like global monopole spacetime in terms of such symmetries is exhibited. \\

\indent  In $2006,$ Hayward \cite{Hayward} modeled the famous exact regular black hole metric, which is the solution of the EFE in spherical symmetry, and it is a simple and singularity free black hole spacetime in general relativity. The line element of HBH spacetime, in spherical coordinates $(t, r, \theta, \phi)$,  is given by
\bea\label{HM}
ds^2=-\left( 1-\frac{2mr^2}{r^3+2mb^2}\right)dt^2 +\left( 1-\frac{2mr^2}{r^3+2mb^2}\right)^{-1} dr^2  + r^2\left( d\theta^2 + \sin^2\theta d\phi^2\right),  
\eea
where the parameters $m$ and $b$ represents mass and length-scale respectively. The metric (\ref{HM}) is non-singular, because if $r \to \infty,$ the metric approaches to $1-\frac{2m}{r}$, and if $r\to 0$, it approaches unity smoothly. The metric  consists of the least number of free parameters ($b$ only) with the properties
\begin{enumerate}
	\item Schwarzschild  asymptotic behavior at large radii and
	
	\item regularly at the center such that $F(r)=1-\frac{2mr^2}{r^3+2mb^2} \to 1+O(r^2),$.
\end{enumerate}
Hence it is minimal. The importance of Hayward spacetime \eqref{HM} is realized from the several studies, such as,  Chiba and Kimura \cite{CK} have obtained the  timelike geodesics and null geodesics equation of a particle in HBH spacetime and   Maluf and Neves \cite{MN} have studied certain thermodynamic quantities like  Hawking temperature, entropy and heat capacity of HBH spacetime. The stability of the thin-shell wormholes constructed by the HBH spacetime are studied by Halilsoy et al. \cite{HOM}.  However, several curvature properties of HBH spacetime are yet to be investigated. \\
%\indent Recently, in $2020,$ Shaikh et al. \cite{SDHK_interior_2020} studied the curvature properties of interior black hole metric. It is found that an interior black hole metric admits several non-trivial geometric structures, such as, pseudosymmetry, pseudosymmetry due to conformal curvature, pseudosymmetry due to concircular curvature, pseudosymmetry due to conharmonic curvature, Ricci pseudosymmetry due to projective curvature and Roter type etc., which are also obtained in HBH spacetime. 
%================================================================================
%\indent However, we are with incomplete knowledge about the curvature properties of Haywards spacetime. Thats why the present study is devoted to deduce the geometric structures of Haywards black hole metric in terms of curvature restrictions. \\

\indent  The  purpose of the article is to focus and to determine  of several geometric properties of HBH spacetime.  %It is noteworthy to mention that the spacetime \eqref{HM}  realizes various  important of geometric structures in the differential geometry such as,  it is a Roter type manifold and conformal curvature $2$-forms are recurrent. Also, it is neither quasi-Einstein nor Einstein but it is a $2$-quasi-Einstein  and Einstein spacetime of level $2$ and Ricci tensor is compatible for $C$, $P$, $R$, $W$ and $K$. 
It is found that the HBH spacetime is  neither semisymmetric nor Ricci generalized pseudosymmetric, but it is pseudosymmetric and satisfies several pseudosymmetric type curvature conditions, such as, pseudosymmetry due to  concircular, conharmonic and conformal curvature tensors. Also, we have exhibited that  both $Q(g,C)$ and $Q(S,C)$ are linearly dependent on the difference $(C\cdot R-R\cdot C)$. It is also proved that the HBH spacetime is an Einstein manifold of level $2$, $2$-quasi Einstein, generalized quasi-Einstein and  Roter type manifold. Moreover, the nature of the stress energy momentum tensor of HBH spacetime is pseudosymmetric.\\

\indent The article is embellished as follows: we discuss some definitions of geometric structures in in Section $2$, which are essential throughout the paper to investigate the geometric properties of HBH spacetime. Section $3$ deals with the study of the HBH spacetime and obtained some interesting results. In Section $4$, certain geometric properties of energy momentum tensor of HBH spacetime are determined. Based on the curvature properties of HBH spacetime, an worthy comparison with the point-like global monopole spacetime  and Reissner-Nordstr\"om spacetime has been exhibited in Section $5$. Section $6$ is concerned with the nature of Ricci soliton and Ricci-Yamabe soliton admited by the HBH spacetime. Section $7$ is devoted to a comparative study between between the HBH spacetime and the point-like global monopole spacetime with respect to different kinds of symmetry, such as, motion, curvature collineation and curvature inheritance.

\section{\bf Preliminaries}
%%%%%%%%%%%%%%%%%%%%%%%%%%%%%%%%%%%%%%%%%%%%%%%%%%%%%%%
%%%%%%%%%%%%%%%%%%%%%%%%%%%%%%%%%%%%%%%%%%%%%%%%%%%%%%%
This section consists of various rudimentary facts about  various geometric structures, Ricci soliton and symmetries (such as,  motion, curvature collineation (also, curvature inheritance) for the (1,3)-type curvature tensor and for the (0,4)-type curvature tensor, Ricci collineation and Ricci inheritance), which are necessary for investigating the geometric structures on HBH spacetime. \\

The Kulkarni-Nomizu product  $A \wedge U$ of two (0,2)-type symmetric tensors $A$ and $U$ is defined as (\cite{DGHS11, Glog02, G08, Kowa06}):
$$(A\wedge U)_{pq\mu\nu}=A_{p\nu}U_{q\mu}-A_{p\mu}U_{q\nu}+A_{q\mu}U_{p\nu}-A_{q\nu}U_{p\mu}.$$

For $j=1,2,3,4$ we will consider $\varpi,\varpi_j   \in \chi(M),$ Lie algebra of all smooth vector fields throughout the paper.  For a symmetric $(0,2)$-type tensor $Z$, the endomorphism $h_1\wedge_Z h_2$   can be defined as (\cite{DDHKS00, DD91, DGHS11})
$$(\varpi_1\wedge_Z \varpi_2)v=Z(\varpi_2,v)\varpi_1-Z(\varpi_1,v)\varpi_2.$$
Now, the endomorphisms $\mathcal{R}, \mathcal{W}, \mathcal{C},\mathcal{K}$ and $\mathcal{P}$ can be defined on $M$ as follows: (\cite{DHJKS14,SC21,SK18,SK19})
$$\mathcal{R}(\varpi_1,\varpi_2)=[\nabla_{\varpi_1},\nabla_{\varpi_2}]-\nabla_{[\varpi_1,\varpi_2]},$$
$$\mathcal{W}(\varpi_1,\varpi_2)=\mathcal{R}(\varpi_1,\varpi_2)-\frac{\kappa}{n(n-1)}\varpi_1\wedge_g \varpi_2,$$
$$\mathcal{C}(\varpi_1,\varpi_2)=\mathcal{R}(\varpi_1,\varpi_2)-\frac{1}{n-2}\left(\varpi_1\wedge_g \mathcal{J}\varpi_2 +\mathcal{J}\varpi_1 \wedge_g \varpi_2-\frac{\kappa}{n-1}\varpi_1\wedge_g \varpi_2\right),$$
$$\mathcal{K}(\varpi_1,\varpi_2)=\mathcal{R}(\varpi_1,\varpi_2)-\frac{1}{n-2}\left(\varpi_1\wedge_g \mathcal{J} \varpi_2+\mathcal{J}\varpi_1\wedge_g \varpi_2\right),$$
$$\mathcal{P}(\varpi_1,\varpi_2)=\mathcal{R}(\varpi_1,\varpi_2)-\frac{1}{n-1}\varpi_1\wedge_S \varpi_2,$$
where $\mathcal{J}$ is the Ricci operator defined as $S(\varpi_1,\varpi_2)=g(\varpi_1,\mathcal{J}\varpi_2)$.  Now, we define the $(0,4)$-type tensor field $T$ corresponding to the endomorphism $\mathscr{T}(\varpi_1,\varpi_2)$ on $M$ as 
$$T(\varpi_1,\varpi_2,\varpi_3,\varpi_4)=g(\mathscr{T}(\varpi_1,\varpi_2)\varpi_3,\varpi_4).$$
If the endomorphism $\mathscr{T}$ is replaced  by  $\mathcal{R}$ (resp., $\mathcal{W}$, $\mathcal{K}$,  $\mathcal{C}$ and $\mathcal{P}$) in above, we obtain the Riemann (resp., concircular, conharmonic, conformal  and projective) curvature  tensor $R$ (resp.,  $W$,  $K$,  $C$ and  $P$)  of type $(0,4)$.
These tensors  are locally given by
$$R_{pq\mu \nu}=g_{p\alpha}(\partial_n \Gamma^\alpha_{q\mu}-\partial_\mu \Gamma^\alpha_{q\nu}+\Gamma^\beta_{q\mu}\Gamma^\alpha_{\beta \nu}-\Gamma^\beta_{q\nu}\Gamma^\alpha_{\beta \mu}),$$
$$C_{pq\mu \nu}=R_{pq\mu\nu}+\frac{\kappa}{2(n-1)(n-2)} (g\wedge g)_{pq\mu\nu}-\frac{1}{n-2}(g\wedge S)_{pq\mu\nu},$$
$$W_{pq\mu\nu}=R_{pq\mu\nu}-\frac{\kappa}{2n(n-1)}(g\wedge g)_{pq\mu\nu},$$
$$K_{pq\mu\nu}=R_{pq\mu\nu}-\frac{1}{n-2}(g\wedge S)_{pq\mu\nu} \ and \ $$ 
$$P_{pq\mu\nu}=R_{pq\mu\nu}-\frac{1}{n-1}(g_{p\nu}S_{q\mu}-g_{q\nu}S_{p\mu}),$$
where $\partial_\alpha=\frac{\partial}{\partial x^\alpha}$ and $\Gamma^\alpha_{q\mu} $ denotes the  Christoffel symbols of $2$nd kind.

Let $H$ be a $(0,k)$-type ($k\geq 1$) tensor on  $M$. Then the $(0,k+2)$-type tensor $T\cdot H$  is given by  \cite{DGHS98, DH03, SK14}
$$(T\cdot H)_{q_{_1}q_{_2}\cdots q_{_k}\nu\mu}=-g^{pr}[T_{\mu\nu q_{_1}r}H_{pq_{_2}\cdots q_{_k}}+\cdots +T_{\mu \nu q_{_k} r}H_{q_{_1}q_{_2}\cdots p}].$$ 
Again, for a symmetric $(0,2)$-type tensor field $Z$, the Tachibana tensor $Q(Z,H)$ of type $(0,k+2)$ is obtained  as follows:  \cite{DGPSS11,SDHJK15,Tach74}
$$Q(Z,H)_{q_{_1}q_{_2}\cdots q_{_k} \nu\mu}=Z_{\mu q_{_1}}H_{\nu q_{_2}\cdots q_{_k}}+ \cdots +Z_{\mu q_{_k}}H_{q_{_1}q_{_2}\cdots \nu}-Z_{\nu q_{_1}}H_{\mu q_{_2}\cdots q_{_k}}-  \cdots -Z_{\nu q_{_k}}H_{q_{_1}q_{_2}\cdots \mu}.$$

\begin{defi} \cite{AD83, Cart46,  Desz92, Desz93, DGHZ15, DGHZ16,SAAC20N, SK14, SKppsnw, Szab82, Szab84, Szab85} 
Let $M$ be a semi-Riemannian manifold. $M$ is called a H-semisymmetric type manifold due to T if $M$ possesses the relation $T\cdot H=0$. Further, $M$ is said to be H-pseudosymmetric type manifold due to T if the relation $T\cdot H=\mathcal{L}_H Q(Z,H)$ holds for a smooth function $\mathcal{L}_H$  on $\{ x\in M: Q(Z,H)\neq 0 \ at \ x \}$ (i.e.,  the tensors $T\cdot H$ and $Q(Z,H)$ are linearly dependent).
\end{defi} 
In the above definition, if we replace $T$ = $R$ and $H$ = $R$ (resp., $P$, $K$, $W$, $C$ and $S$), then the $H$-semisymmetric type manifold due to $T$ turns into semisymmetric  (resp., projectively, conharmonically, concircularly, conformally, Ricci semisymmetric) manifold  and if $T=R$, $H=R$ and $Z=g$ (resp., $S$), then the $H$-pseudosymmetric type manifold due to $T$ becomes Deszcz pseudosymmetric (resp., Ricci generalized pseudosymmetric) manifold.  Also, if we replace $T$ = $W$, $C$, $P$ and $K$, then we obtain  several pseudosymmetric type curvature conditions.

\begin{defi}$($\cite{DGHZ16, DGJZ-2016, DGP-TV-2011, S09, SKH11, SK19,SYH09}$)$ 
 $M$ is called quasi-Einstein (resp., Einstein and $2$-quasi-Einstein) manifold if for a scalar $\alpha$ the rank of $(S-\alpha g)$ is $1$ (resp., $0$ and $2$).  In particular, for $\alpha =0$ the quasi-Einstein manifold turns into Ricci simple.
A generalized quasi-Einstein manifold (in the sense of Chaki \cite{C01}) is defined as $$S=\alpha g+\beta \Pi \otimes \Pi +\gamma (\Pi \otimes \phi +\phi \otimes \Pi)$$  where $\alpha$, $\beta$ and $\gamma$ are scalars and $\Pi$, $\phi$ are $1$-forms.
\end{defi}

It may be mentioned that Robertson Walker spacetime \cite{ARS95, Neill83, SKMHH03} is quasi-Einstein, Kaigorodov spacetime \cite{SDKC19} is Einstein, Kantowski-Sachs spacetime \cite{SC21} and Som-Raychaudhuri spacetime \cite{SK16srs}  are $2$-quasi-Einstein, Vaidya metric \cite{SKS19}, G\"odel spacetime \cite{DHJKS14} and Morris-Thorne spacetime \cite{ECS22} are also a Ricci simple manifold. 
\begin{defi}
If the Ricci tensor $S$ of a semi-Riemannain manifold $M$ satisfies the relation 
$$(\nabla_{\varpi_{1}}S)(\varpi_2,\varpi_3) + (\nabla_{\varpi_{2}}S)(\varpi_3,\varpi_1) + (\nabla_{\varpi_{3}}S)(\varpi_1,\varpi_2)=0,$$ 
 $$( resp., (\nabla_{\varpi_{1}}S)(\varpi_2,\varpi_3) = (\nabla_{\varpi_{2}}S)(\varpi_1,\varpi_3))$$
   then it is known as cyclic parallel Ricci tensor  (see, \cite{Gray78, SB08, SJ06, SJ07}) (resp., Codazzi type Ricci tensor (see \cite{F81, S81})). 
\end{defi}
It may be noted that the   Ricci tensor of $(t-z)$-type plane wave spacetime \cite{EC21} is of Codazzi type and the Ricci tensor of cyclic parallel has been found in G\"odel spacetime \cite{DHJKS14}.
\begin{defi} $($\cite{Bess87, SK14, SK19}$)$ 
	A semi-Riemannian manifold $M$ is an Einstein manifold  of level $4$ (resp., $3$ and $2$) if it satisfies  
	 $$\vartheta_1 g+ \vartheta_2 S+ \vartheta_3 S^2+ \vartheta_4 S^3+  S^4=0,$$
%\eea
(resp., $\vartheta_5 g+ \vartheta_6 S+  \vartheta_7 S^2+ S^3=0$ \ and \ $\vartheta_8 g+ \vartheta_9 S+  S^2=0$),
where $\vartheta_i$ $(1\leq i \leq 9)$ are smooth functions on $M$.
\end{defi}

We mention  that Vaidya-Bonner spacetime \cite{SDC} and Lifshitz spacetime \cite{SSC19} are Ein$(3)$ while Siklos spacetime \cite{SDKC19} and Nariai spacetime \cite{SAAC20N} are Ein$(2)$ manifolds.
\begin{defi} 	 
	If the Riemann tensor $R$  can be written in the form  
	\bea
	R = S^2 \wedge (\varsigma_6 S^2)+S \wedge (\varsigma_4S+\varsigma_5 S^2)+g \wedge (\varsigma_1 g+\varsigma_2 S+\varsigma_3 S^2)\nonumber 
	\eea
 for some scalars $\varsigma_i$ $(1\leq i\leq 6)$, then $M$ is called generalized Roter type manifold \cite{Desz03, DGJPZ13, DGJZ-2016, DGP-TV-2015, SK16,SK19}.   Further, $M$ is known as a Roter type  manifold \cite{Desz03, DG02, DGP-TV-2011, DPSch-2013, Glog-2007} if  $g\wedge g,$ $S\wedge S$ and $g\wedge S$  are linearly dependent on $R$ (i.e., $\varsigma_3=\varsigma_5=\varsigma_6=0$).	
%	\[\left(resp.,\;R=g \wedge (\mu_7 g+\mu_8 S+\mu_9 S^2 % \right)\] 
	
%	where $\mu_i (1\leq i\leq 9)$ are some scalars.

\end{defi}
It may be noted that Nariai spacetime \cite{SAAC20N}, Melvin magnetic metric \cite{SAAC20} as well as Robinson-Trautman spacetime \cite{SAA18} are Roter type, while Lifshitz metric \cite{SSC19} and Vaidya-Bonner metric \cite{SDC} are generalized Roter type manifolds.
\begin{defi} \cite{TB89, TB93}
	A weakly $T$-symmetric manifold $M$ is defined by the equation
	\beb
	(\nabla_{\varpi} T)(\varpi_1,\varpi_2,\varpi_3,\varpi_4)&=& \Pi(\varpi)\otimes T(\varpi_1,\varpi_2,\varpi_3,\varpi_4)+ \Omega_1(\varpi_1)\otimes T(\varpi,\varpi_2,\varpi_3,\varpi_4)\\ &+& \Omega_1(\varpi_2)\otimes T(\varpi_1,\varpi,\varpi_3,\varpi_4)+ \Omega_2(\varpi_3)\otimes T(\varpi_1,\varpi_2,\varpi,\varpi_3)\\ &+& \Omega_2(\varpi_4)\otimes T(\varpi_1,\varpi_2,\varpi_3,\varpi),
	\eeb
	where  $\Pi$, $\Omega_1$, $\Omega_2$  are $1$-forms on $M$.
	  In particular, $M$ reduces to a recurrent \cite{Ruse46,Ruse49a,Walk50}  (resp., Chaki pseudosymmetric  \cite{Chak87, Chak88})  manifold if  $\Omega_1$=$\Omega_2$=$0$    (resp., $\Omega_1$=$\Omega_2=\Pi\slash2$).
\end{defi}

\begin{defi}$($\cite{DD91,DGJPZ13, MM12a, MM12b, MM13}$)$ Let $T$ be a $(0,4)$-type tensor field and $\mathcal{Z}$ be the endomorphism corresponding to a tensor $Z$ of type $(0,2)$ on $M$. Then, the tensor $Z$ is said to be $T$-compatible if $M$ admits 
	\beb
	\mathop{\mathcal{S}}_{\varpi_1,\varpi_2,\varpi_3} T(\mathcal{Z}\varpi_1,\varpi,\varpi_2,\varpi_3)=0,
	\eeb
	where the cyclic sum over  $\varpi_1$, $\varpi_2$ and $\varpi_3$ is denoted by $\mathcal{S}$.	Again, $T$-compatibility of an $1$-form $\zeta$ is defied by the $T$-compatibility of $\zeta \otimes \zeta$.\end{defi}

	In the above definition, if we replace $T$ by $R$ (resp., $P$, $K$, $W$ and $C$) then  we obtain Riemann (resp., projective, conharmonic, concircular and conformal) compatibility of $Z$.

\begin{defi}\label{defi2.8} 
Let $T$ be a $(0,4)$-type tensor on $M$ and $\mathcal{Z}$ be the endomorphism corresponding to a $(0,2)$-type tensor $Z$ on $M$. If $M$ possesses the relation   
	\beb
	\mathop{\mathcal{S}}_{\varpi_1,\varpi_2,\varpi_3}(\nabla_{\varpi_1}T)(\varpi_2,\varpi_3,\varpi_4,\varpi_5)=\mathop{\mathcal{S}}_{\varpi_1,\varpi_2,\varpi_3}\Sigma(\varpi_1)T(\varpi_2,\varpi_3,\varpi_4,\varpi_5)
	\eeb
	  for an $1$-form $\Sigma$, then the curvature $2$-forms $\Omega^m_{(T)l}$ \cite{LR89} are recurrent \cite{MS12a,MS13a,MS14}. Further,   the $1$-forms  $\Lambda_{(Z)l}$\cite{SKP03} are recurrent if 
	$$(\nabla_{\varpi_1}Z)(\varpi_2,\varpi)-(\nabla_{\varpi_2}Z)(\varpi_1,\varpi)=\Sigma(\varpi_1)  Z(\varpi_2,\varpi)-\Sigma(\varpi_2) Z(\varpi_1,\varpi)$$ holds on $M$ for an $1$-form $\Sigma$.
\end{defi}
%

%
%\begin{defi}\label{2.6}
%	The Ricci tensor of a semi-Riemannian manifold $M$, which is defined by the relation 
	%\[(\nabla_{W_1}S)(W_2,W_3)=(\nabla_{W_2}S)(W_1,W_3)\]
	% \[\left(resp.,\;  \mathop{\mathcal{S}}_{W_1,W_2,W_3}(\nabla_{W_1}S)(W_2,W_3)=0\right)\] holds, then it is called Codazzi type %\cite{F81, S81}  (resp., cyclic parallel \cite{ Gray78, SB08, SS06, SS07,SSC19}) Ricci tensor, 
	 % where $\mathcal{S}$ is the cyclic sum over the vector fields $W_1, W_2$ and $W_3$.
%\end{defi}
%We mention that the $(t-z)$-type plane wave metric \cite{EC21} has Ricci tensor of Codazzi type  and G\"odel spacetime \cite{DHJKS14} is equipped with cyclic parallel Ricci tensor. 
%
%

%

%
\begin{defi}$($\cite{P95, Venz85}$)$
	Let T be a (0,4)-type tensor on $M$ and $L(M)$ be the set of all 1-forms $\Pi$ on $M$ satisfying 	\beb
	\mathop{\mathcal{S}}_{\varpi_1,\varpi_2,\varpi_3}\Pi(\varpi_1)\otimes T(\varpi_2,\varpi_3,\varpi_4,\varpi_5)=0
	\eeb
	with $dimL(M)\geq 1$. Then $M$ is called a $ T $-space by Venzi. 
	\end{defi}
%
%%%%%%%%%%%%%%%%%%%%%%%%%%%%%%%%%%%%%%%%%%%%%%%%%%%%%%%%%%%%%%%%%%%%%%%%%%

Several notions of geometrical symmetries, such as, motion, curvature collineation for (0,4)-type curvature tensor and for (1,3)-type curvature tensor, curvature inheritance for (0,4)-type curvature tensor and for (1,3)-type curvature tensor, Ricci collineation and Ricci inheritance, all of which are originated from the Lie derivatives of different tensors, are essential to be reviewed for the study of symmetry in the HBH spacetime.

\begin{defi}
	A manifold $M$ admits motion with respect to some vector field $\xi$ if $\pounds_\xi g=0$. The vector field $\xi$ is also called Killing.
\end{defi}

 Katzin {et al.} \cite{KLD1969,KLD1970}, in 1969,  introduced the concept of curvature collineation for the (1,3)-type curvature tensor by vanishing Lie derivative of the (1,3)-type Riemann curvature tensor with respect to some vector field. Again, in 1992, by introducing the notion of curvature inheritance for the (1,3)-type curvature tensor, Duggal \cite{Duggal1992} generalizes the concept of curvature collineation.

\begin{defi}\label{def_CI} (\cite{Duggal1992})
	A  semi-Riemannian manifold $M$ admits curvature inheritance for the (1,3)-type curvature tensor $\widetilde R$ if $M$ satisfies 
	$$\pounds_\xi \widetilde R=\lambda \widetilde R$$
	for a non-Killing vector field $\xi$, 
	where $\lambda$ is a scalar function and the (1,3)-type curvature tensor $\widetilde R$ is related to the (0,4)-type curvature tensor $R$  by   $R(v_1,v_2,v_3,v_4)=g(\widetilde R(v_1,v_2)v_3,v_4)$. In particular, if $\lambda=0$, then it turns into curvature collineation \cite{KLD1969,KLD1970} for the (1,3)-type curvature tensor $\widetilde R$ (i.e., $\pounds_\xi \widetilde R=0$). 
\end{defi}

 \begin{defi} \label{def_RI}(\cite{Duggal1992})
	A semi-Riemannian manifold $M$ realizes Ricci inheritance if for some vector field $\xi$ and for some scalar function $\lambda$, $M$ possesses the relation
	$$\pounds_\xi S=\lambda S.$$
	. Further, if $\lambda=0$, it transforms into Ricci collineation (i.e., $\pounds_\xi S=0$).
\end{defi}

Recently,  generalizing the notion of curvature inheritance for (0,4)-type curvature tensor $R$ Shaikh and Datta \cite{ShaikhDatta2022} introduced the concept of generalized curvature inheritance for (0,4)-type curvature tensor $R$, which is given as follows:
\begin{defi}\label{def_GCI} (\cite{ShaikhDatta2022})
	A semi-Riemannian manifold $M$ admits generalized curvature inheritance for (0,4)-type curvature tensor $R$ if there is a non-Killing vector field $\xi$ which satisfies the relation
	$$\pounds_\xi  R=\lambda  R + \lambda_1 g\wedge g +\lambda_2 g\wedge S +\lambda_3 S\wedge S ,$$
	where $\lambda,\lambda_1,\lambda_2,\lambda_3$ are the scalar functions. In particular, if $\lambda_i=0$ for $i=1,2,3$, then  $M$ admits curvature inheritance  for (0,4)-type curvature tensor $R$. Further, if  $\lambda=0=\lambda_i$ for $i=1,2,3$, then it becomes curvature collineation  for (0,4)-type curvature tensor $R$.
\end{defi}

\section{\bf Hayward black hole spacetime admitting geometric structures} 
%%%%%%%%%%%%%%%%%%%%%%%%%%%%%%%%%%%%%%%%%%%%%%%%%%%%%%%%%%%%%%%%%%%%%%%
%=======================================================================
%
In coordinates  $(t,r,\theta,\phi)$, the metric tensor  of HBH spacetime is given by:
$$
g=\left(
\begin{array}{cccc}
	-(1- \frac{2mr^2}{r^3+2mb^2}) & 0 & 0 & 0 \\
	0 & (1- \frac{2mr^2}{r^3+2mb^2})^{-1} & 0 & 0 \\
	0 & 0 & r^2 & 0 \\
	0 & 0 & 0 & r^2 \sin^2\theta
\end{array}
\right).
$$
Now, the components of the metric $g$ are  
$$\begin{array}  {c}g_{11}=-\left(1-\frac{2mr^2}{r^3+2mb^2}\right),\; g_{22}=\left(1-\frac{2mr^2}{r^3+2mb^2}\right)^{-1}, \\ g_{33}=r^2, \;g_{44}=r^2\sin^2\theta, \; g_{ij}=0, \;\mbox{otherwise}.
\end{array}$$

Let $B=2b^2m+r^2(r-2m)$, $B_1=2b^2m+r^3$, $B_2=4b^2m-r^3$, $B_3=b^2m-r^3$ and $B_4=10b^2m-r^3$. The non-vanishing components of the  Christoffel symbols $(\Gamma^h_{ij})$ of 2nd kind are calculated as  given  below:
$$\begin{array}{c}
	\Gamma^2_{11}=-\frac{mrBB_2}{B_1^3}, \ \ 
	\Gamma^1_{12}=-\frac{mrB_2}{BB_1}=-\Gamma^2_{22}, \\
	\Gamma^3_{23}=\frac{1}{r}=\Gamma^4_{24}, \ \ 
	\Gamma^2_{33}=-r+\frac{2mr^3}{B_1}, \\
	\Gamma^4_{34}= \cot\theta, \ \ 
	\Gamma^2_{44}=-\frac{rB\sin^2\theta}{B_1}, \\
	\Gamma^3_{44}=-\cos\theta \sin\theta. 
\end{array}$$

The non-vanishing components  of the  Riemann-curvature $(R_{abcd})$ and Ricci tensor $(S_{ab})$ and the scalar curvature $\kappa$ are calculated as  given  below:
$$\begin{array}{c}
	R_{1212}=-\frac{2m(2b^2m(2b^2m-7r^3)+r^6)}{B_1^3}, \ \ 
	
	R_{1313}=\frac{mr^2BB_2}{B_1^3}=\frac{1}{\sin^2\theta}R_{1414}, \\
	
	R_{2323}=\frac{mr^2B_2}{BB_1}=\frac{1}{\sin^2\theta}R_{2424}, \  \
	R_{3434}=\frac{2mr^4 \sin^2\theta}{B_1}; \\
\end{array}$$
$$\begin{array}{c}
	S_{11}=\frac{24b^2m^2BB_3}{B_1^4}, \ \ 
	S_{22}=\frac{24b^2m^2(-b^2m+r^3)}{BB_1^2}, \\
	S_{33}=-\frac{12b^2m^2r^2}{B_1^2}, \ \ 
	S_{44}=\sin^2\theta  S_{33};
\end{array}$$
 $$\kappa=\frac{24b^2m^2(r^3-4b^2m)}{(2b^2m+r^3)^3}.$$ 
 From the above calculation, one can obtain the following:
\begin{pr} \label{pr1}
	The HBH spacetime is neither Einstein nor quasi-Einstein but $(i)$ it is $2$-quasi-Einstein %as rank $(S-\alpha g)=2$
	 for $\alpha=-\frac{12b^2m^2}{B_1^2 } $ and $(ii)$  for $\alpha=-\frac{12b^2m^2}{B_1^2}$,  $\beta=1$, $\gamma=1$, $\Pi$=$\left\lbrace -\frac{B}{B_1},1,0,0\right\rbrace $ and $\phi$= $\left\lbrace \frac{36b^2m^2r^3+B_1^2B}{2B_1^3},\frac{18b^2m^2r^3}{B_1^2B }-\frac{1}{2},0,0\right\rbrace $,    it is generalized quasi-Einstein in the sense of Chaki.
\end{pr}
Let $\mathcal{K}^1=(g\wedge g)$, $\mathcal{K}^2=(g\wedge S)$ and $\mathcal{K}^3=(S\wedge S).$ Then the components other than zero of   $\mathcal{K}^1$, $\mathcal{K}^2$ and $\mathcal{K}^3$ are calculated as  given  below:
$$\begin{array}{c}
	\mathcal{K}^1_{1212}=2, 
	\mathcal{K}^1_{1313}=\frac{2r^2B}{B_1}=\frac{1}{\sin^2\theta}\mathcal{K}^1_{1414}, \\
	\mathcal{K}^1_{2323}=-\frac{2r^2B_1}{B}=\frac{1}{\sin^2\theta}\mathcal{K}^1_{2424}, \ \
	\mathcal{K}^1_{3434}=-2r^4 \sin^2\theta; 
\end{array}$$

$$\begin{array}{c}
	\mathcal{K}^2_{1212}=-\frac{48b^2m^2B_3}{B_1^3}, \ \
	\mathcal{K}^2_{1313}=-\frac{12b^2m^2r^2BB_2}{B_1^4}=\frac{1}{\sin^2\theta}\mathcal{K}^2_{1414}, \\
	
	\mathcal{K}^2_{2323}=\frac{12b^2m^2r^2B_2}{BB_1^2}=\frac{1}{\sin^2\theta}\mathcal{K}^2_{2424}, \ \
	\mathcal{K}^2_{3434}=\frac{24b^2m^2r^4\sin^2\theta}{B_1^2};
\end{array}$$

$$\begin{array}{c}	
	\mathcal{K}^3_{1212}=\frac{1152b^4m^4B_3^2}{B_1^6}, \ \ 
	\mathcal{K}^3_{1313}=\frac{576b^4m^4r^2B_3B}{B_1^6}=\frac{1}{\sin^2\theta}\mathcal{K}^3_{1414}, \\
	
	\mathcal{K}^3_{2323}=-\frac{576b^4m^4r^2B_3}{B_1^4B}=\frac{1}{\sin^2\theta}\mathcal{K}^3_{2424}, \ \ 
	\mathcal{K}^3_{3434}=-\frac{288b^4m^4r^4\sin^2\theta}{B_1^4}. \\
\end{array}$$
\indent From the above calculation, it follows that  $S\wedge S$, $g\wedge S$, $g\wedge g$,  and   $R$ are linearly dependent in HBH spacetime, and hence
the Riemann tensor $R$ can be explicitly given as follows:
\bea\label{RT}
R=\varsigma_1 \mathcal{K}^1 + \varsigma_2 \mathcal{K}^2 + \varsigma_3 \mathcal{K}^3
\eea
where $\varsigma_1=m(\frac{2}{3r^3}-\frac{1}{B_1})$, $\varsigma_2=\frac{1}{36}(10+\frac{16b^2m}{r^3}+\frac{r^3}{b^2m})$ and $\varsigma_3=\frac{B_2B_1^3}{432b^4m^3r^3}$.\\  On contraction the relation  \eqref{RT} entails 
\bea\label{Ein2}
S^2 + \vartheta_1 S + \vartheta_2 g=0
\eea
where $\vartheta_1=\frac{12b^2m^2B_2}{B_1^3}$ and $\vartheta_2=\frac{288b^4m^4B_3}{B_1^5}$. \\
\indent From the  relation \eqref{RT} and \eqref{Ein2}, we can state the following: 
\begin{pr}\label{pr2}
	The HBH spacetime is neither Ein$(3)$ nor generalized Roter type   but it fulfills $(i)$  Roter type     and $(ii)$ Einstein manifold of level $2.$   
\end{pr}
The non-vanishing components  $C_{abcd}$ of the conformal curvature tensor $C$ (upto symmetry) are calculated and given as below:
$$\begin{array}{c}
	C_{1212}=\frac{2mr^3B_2}{B_1^3}, \ \ 
	C_{1313}=-\frac{mr^5B_2B}{B_1^4}=\frac{1}{\sin^2\theta}C_{1414}, \\
	
	C_{2323}=\frac{mr^5B_2}{B_1^2B}=\frac{1}{\sin^2\theta}C_{2424}, \ \ 
	C_{3434}=-\frac{2mr^7B_2\sin^2\theta}{B_1^3}.\\
\end{array}$$
%\indent Calculated the non-zero components of the covariant derivatives of the tensor $R_{abcd}$ and $C_{abcd}$ given as below: \\
%Let $T^1_{abcd,f}= \nabla_f R_{abcd}$ and  $T^2_{abcd,f}= \nabla_f C_{abcd}$
If  $\mathscr D_{abcd,f}= \nabla_f R_{abcd}$ and $\mathscr F_{abcd,f}= \nabla_f C_{abcd}$, 
then the components other than zero of $\nabla{R}$ and $\nabla{C}$  are obtained  as follows:
$$\begin{array}{c}
	\mathscr D_{1212,2}= \frac{6m(40b^4m^2r^2-32b^2mr^5+r^8)}{B_1^4}, \ \ 
	\mathscr D_{1213,3}=\frac{3mr^4B_4B}{B_1^4}=\mathscr D_{1313,2}, \\
	\mathscr D_{1214,4}=-\frac{3mr^4B_4B\sin^2\theta}{B_1^4}=\mathscr D_{1414,2}, \ \ 
	
	\mathscr D_{2323,2}=-\frac{3mr^4B_4}{B_1^2B}=\frac{1}{\sin^2\theta}	\mathscr D_{2424,2}, \\
	
	\mathscr D_{2334,4}=\frac{3mr^6 \sin^2\theta}{B_1^2}=-\mathscr D_{2434,3}=-\frac{1}{2} \mathscr D_{3434,2}; 
\end{array}$$

$$\begin{array}{c}	
	\mathscr F_{1212,2}=\frac{6mr^2(8b^4m^2-12b^2mr^3+r^6)}{B_1^4}=-\frac{1}{r^4 \sin^2\theta}\mathscr F_{3434,2}, \\
	\mathscr F_{1213,3}=-\frac{3mr^4BB_2}{B_1^4}=\frac{1}{\sin^2\theta}\mathscr F_{1214,4}, \ \
	
	\mathscr F_{1313,2}=-\frac{3mr^4(8b^4m^2-12b^2mr^3+r^6)B}{B_1^5}=\frac{1}{\sin^2\theta}\mathscr F_{1414,2}, \\
	
	\mathscr F_{2323,2}=\frac{3mr^4(8b^4m^2-12b^2mr^3+r^6)}{B_1^3B}=\frac{1}{\sin^2\theta}\mathscr F_{2424,2}, \ \ 
	\mathscr F_{2334,4}=-\frac{3mr^6B_2\sin^2\theta}{B_1^3}=-\mathscr F_{2434,3}. 
\end{array}$$
From the above components we get the following proposition:
\begin{pr}\label{pr3}
	The HBH spacetime is not conformally recurrent but its $(i)$  conformal $2$-form are recurrent for the $1$-forms  $\{ 0, -\frac{6b^2m(8b^2m-5r^3)}{8b^4m^2r+2b^2mr^4-r^7},0 ,0 \}$ and $(ii)$ the general form of  $R$-compatible tensor and $C$-compatible tensor are given by
	
	$$
	\left(
	\begin{array}{cccc}
		\mathscr{Z}_{11} &\mathscr{Z}_{12} & 0 & 0 \\
		\mathscr{Z}_{12} & \mathscr{Z}_{22} & 0 & 0 \\
		0 & 0 & \mathscr{Z}_{33} & \mathscr{Z}_{34} \\
		0 & 0 & \mathscr{Z}_{34} & \mathscr{Z}_{44}
	\end{array}
	\right)
	$$
	where $\mathscr{Z}_{ij}$ are arbitrary scalars.

\end{pr}
Let $\mathcal{M}^{1}=R\cdot R$, $\mathcal{M}^{2}=R\cdot C$, $\mathcal{M}^{3}=C\cdot R$,  $\mathscr P^1=Q(g,R)$, $\mathscr P^2=Q(S,R)$, $\mathscr P^3=Q(g,C)$ and $\mathscr P^4=Q(S,C)$. Then the  components  of  $\mathcal{M}^{1}$, $\mathcal{M}^{2}$, $\mathcal{M}^{3}$, $\mathscr P^1$, $\mathscr P^2$, $\mathscr P^3$ and  $\mathscr P^4$, which do not vanish, are  given upto symmetry as follows:
$$\begin{array}{c}
	\mathcal{M}^1_{1223,13}=-\frac{3m^2r^5(40b^4m^2-14b^2mr^3+r^6)}{B_1^5}=-\mathcal{M}^1_{1213,23}, \ \ 
	\mathcal{M}^1_{1434,13}=-\frac{3m^2r^7BB_2\sin^2\theta}{B_1^5}=-\mathcal{M}^1_{1334,14}, \\
	
	\mathcal{M}^1_{1224,14}=-\frac{3m^2r^5(40b^4m^2-14b^2mr^3+r^6)\sin^2\theta}{B_1^5}=-\mathcal{M}^1_{1214,24}, \ \ 
	\mathcal{M}^1_{2434,23}=-\frac{3m^2r^7B_2\sin^2\theta}{B_1^3B}=-\mathcal{M}^1_{2334,24}; 
\end{array}$$

	$$\begin{array}{c}
	\mathcal{M}^2_{1223,13}=-\frac{3m^2r^5B_2^2}{B_1^5}= -\mathcal{M}^2_{1213,23}, \ \ 
	\mathcal{M}^2_{1434,13}=\frac{3m^2r^7BB_2\sin^2\theta}{B_1^6}=-\mathcal{M}^2_{1334,14}, \\
	
	\mathcal{M}^2_{1224,14}=-\frac{3m^2r^5B_2^2\sin^2\theta}{B_1^5}=- \mathcal{M}^2_{1214,24}, \ \ 
	
	\mathcal{M}^2_{2434,23}=\frac{3m^2r^7B_2^2\sin^2\theta}{B_1^4B}=-\mathcal{M}^2_{2334,24}; 
	
\end{array}$$

$$\begin{array}{c}
	\mathcal{M}^3_{1223,13}=-\frac{3m^2r^8(40b^4m^2-14b^2mr^3+r^6)}{B_1^6}=-\mathcal{M}^3_{1213,23}, \\
	\mathcal{M}^3_{1434,13}=-\frac{3m^2r^{10}BB_2\sin^2\theta}{B_1^6}=-\mathcal{M}^3_{1334,14}, \\
	
	\mathcal{M}^3_{1224,14}=-\frac{3m^2r^8(40b^4m^2-14b^2mr^3+r^6)\sin^2\theta}{B_1^6}=-\mathcal{M}^3_{1214,24}, \\
	\mathcal{M}^3_{2434,23}=-\frac{3m^2r^{10}B_2\sin^2\theta}{B_1^4B}=-\mathcal{M}^3_{2334,24}; 
	\end{array}$$

$$\begin{array}{c}	
	\mathscr P^1_{1223,13}=\frac{3mr^5B_4}{B_1^3}=-\mathscr P^1_{1213,23}, \ \ 
	\mathscr P^1_{1434,13}=\frac{3mr^7B\sin^2\theta}{B_1^3}=-\mathscr P^1_{1334,14}, \\
	\mathscr P^1_{1224,14}=\frac{3mr^5B_4\sin^2\theta}{B_1^3}=-\mathscr P^1_{1214,24}, \ \ 
	\mathscr P^1_{2434,23}=\frac{3mr^7 \sin^2\theta}{B_1B}=-\mathscr P^1_{2334,24}; \\
\end{array}$$

$$\begin{array}{c}	
	\mathscr P^2_{1223,13}=-\frac{216b^4m^4r^5}{B_1^5}=-\mathscr P^2_{1213,23}, \ \ 
	\mathscr P^2_{1434,13}=-\frac{36b^2m^3r^7B\sin^2\theta}{B_1^5}=-\mathscr P^2_{1334,14}, \\
	\mathscr P^2_{1224,14}=-\frac{216b^4m^4r^5\sin^2\theta}{B_1^5}=-\mathscr P^2_{1214,24}, \ \ 
	\mathscr P^2_{2434,23}=\frac{36b^2m^3r^7\sin^2\theta}{B_1^3B}=-\mathscr P^2_{2334,24}; 
\end{array}$$

$$\begin{array}{c}	
	\mathscr P^3_{1223,13}=\frac{3mr^5B_2}{B_1^3}=-\mathscr P^3_{1213,23}, \ \ 
	\mathscr P^3_{1434,13}=-\frac{3mr^7BB_2\sin^2\theta}{B_1^4}=-\mathscr P^3_{1334,14}, \\
	\mathscr P^3_{1224,14}=\frac{3mr^5B_2\sin^2\theta}{B_1^3}=-\mathscr P^3_{1214,24}, \ \ 
	\mathscr P^3_{2434,23}=-\frac{3mr^7B_2\sin^2\theta}{B_1^2B}=-\mathscr P^3_{2334,24}; 
\end{array}$$

$$\begin{array}{c}	
	\mathscr P^4_{1223,13}=-\frac{72b^4m^4r^5B_2}{B_1^6}= - \mathscr P^4_{1213,23}, \\
	\mathscr P^4_{1434,13}=-\frac{36b^2m^3r^7(2b^2m-r^3)B_2B\sin^2\theta}{B_1^7}=-\mathscr P^4_{1334,14}, \\
	\mathscr P^4_{1224,14}=-\frac{72b^4m^4r^5B_2\sin^2\theta}{B_1^6}=- \mathscr P^4_{1214,24}, \\
	\mathscr P^4_{2434,23}=\frac{36b^2m^3r^7(8b^4m^2-6b^2mr^3+r^6)\sin^2\theta}{B_1^5 B}=-\mathscr P^4_{2334,24}. 
\end{array}$$
From the above components we get the following:
\begin{pr}\label{pr4}
	The HBH spacetime is not Ricci generalized pseodosymmetric but  it is pseudosymmetric and realizes several  pseudosymmetric type curvature relations:
	\begin{enumerate}[label=(\roman*)]
		\item $R\cdot R=-\frac{mB_2}{B_1^2}Q(g,R)$ and hence
		$R\cdot C=-\frac{mB_2}{B_1^2} Q(g,C)$, 
		\item $C\cdot R=-\frac{mr^3B_2}{B_1^3}Q(g,R)$ and hence $C\cdot C=-\frac{mr^3B_2}{B_1^3}Q(g,C),$ 
		\item $R\cdot R-\frac{m(16b^2m-r^3)}{B_2B_1}Q(g,C)=Q(S,R),$
		\item $C\cdot R -R\cdot C=\mathcal{L}_1\ Q(g,R) + \mathcal{L}_2\ Q(S,R),$
		where $\mathcal{L}_1=-\frac{8b^2m^2B_2}{(16b^2m-r^3)B_1^2}$ \ and\ $\mathcal{L}_2=\frac{16b^4m^2-8b^2mr^3+r^6}{(r^3-16b^2m)B_1},$
		\item $C\cdot R -R\cdot C=\mathcal{L}_3\ Q(g,C) + \mathcal{L}_4\ Q(S,C),$
		where $\mathcal{L}_3=\frac{8b^2m^2B_2}{B_1^3}$ and $\mathcal{L}_4=1$.
	\end{enumerate}	
\end{pr}

%\begin{cor}\label{cr1}
%	By Theorem 6.7 of \cite{DGHS11}, since proposition \eqref{pr} imply the following geometric structure :
%	\begin{enumerate}[label=(\roman*)]
%	\item $R\cdot R= \delta_R Q(g,R)$,\ \ $\delta_R=-\frac{m(4b^2m-r^3)}{(2b^2m+r^3)^2},$
%-------------------------------------------------------
%\item $R\cdot C=\delta_R Q(g,C)$, \ $\delta_C=-\frac{m(4b^2m-r^3)}{(2b^2m+r^3)^2},$
%-------------------------------------------------------
%\item $C\cdot R=\delta_C Q(g,R)$, \ \ $\delta_C=\delta_R - (\frac{\kappa}{3}+\frac{\delta_2}{2\delta_3})+\frac{1}{4\delta_3}=-\frac{mr^3(4b^2m-r^3)}{(2b^2m+r^3)^3}$,
%-----------------------------------------------------
%\item $C\cdot C=\delta_CQ(g,C)$, \ $\delta_C=-\frac{mr^3(4b^2m-r^3)}{(2b^2m+r^3)^3}$,
%-----------------------------------------------------
%\item $R\cdot R=Q(S,R)+\delta Q(g,C)$, \ \ $\delta=\delta_R+\frac{\delta_2}{2\delta_3}=\frac{m(16b^2m-r^3)}{(-4b^2m+r^3)(2b^2m+r^3)}$.
%\end{enumerate}
%\end{cor}

%\begin{pr}\label{pr3}
%	Hayward black hole metric \eqref {HM1}, satisfies the pseudosymmetric type curvature relations: 
%	\bea
%	C\cdot R -R\cdot C=\bar{\delta}_1\ Q(g,R) + \bar{\delta}_2\ Q(S,R)
%	\eea
%	where $\bar{\delta}_1=-\frac{8b^2m^2(4b^2m-r^3)}{(16b^2m-r^3)(2b^2m+r^3)^2}$ \ and\ $\bar{\delta}_2=\frac{-16b^4m^2-8b^2mr^3+r^6}{(16b^2m-r^3)(2b^2m+r^3)}$
%
%	\bea
%	\text{and} \ \ \ \ \ C\cdot R -R\cdot C=\bar{\delta}_3\ Q(g,C) + Q(S,C)
%	\eea
%	where $\bar{\delta}_3=\frac{8b^2m^2(4b^2m-r^3)}{(2b^2m+r^3)^3}$.
%\end{pr}
\indent The non-vanishing components  $P_{abcd}$ of the projective curvature tensor $P$ (upto symmetry) of the HBH spacetime are calculated as follows:  
$$\begin{array}{c}
	P_{1212}=\frac{2mr^3B_4}{B_1^3}=-P_{1221}, \ \ 
	P_{1313}=-\frac{mr^5B_4B}{B_1^4}=\frac{1}{\sin^2\theta}P_{1414}, \\
	P_{1331}=-\frac{mr^5B}{B_1^3}=\frac{1}{\sin^2\theta}P_{1441}, \ \ 
	
	P_{2323}=\frac{mr^5B_4}{B_1^2B}=\frac{1}{\sin^2\theta}P_{2424}, \\
	P_{2332}=\frac{mr^5}{B_1B}=\frac{1}{\sin^2\theta}P_{2442}, \ \ 
	P_{3434}=\frac{2mr^7\sin^2\theta}{B_1^2}=-P_{3443}.
\end{array}$$

Let $\mathcal{M}^4=P\cdot S$ and $\mathscr P^5=Q(g,S)$. Then the non-vanishing components of the tensor $\mathcal{M}^4$ and $\mathscr P^5$ are obtained as follows: 
$$\begin{array}{c}
	\mathcal{M}^4_{13,13}=-\frac{36b^2m^3r^5B_2B}{B_1^6}=-\mathcal{M}^4_{13,31},\ \ 
	\mathcal{M}^4_{14,14}=-\frac{36b^2m^3r^5B_2B \sin^2\theta}{B_1^6}=-\mathcal{M}^4_{14,41},\\
	\mathcal{M}^4_{23,23}=\frac{36b^2m^3r^5B_2}{BB_1^4}=-\mathcal{M}^4_{23,32},\ \ 
	\mathcal{M}^4_{24,24}=\frac{36b^2m^3r^5B_2\sin^2\theta}{BB_1^4}=-\mathcal{M}^4_{24,42};
\end{array}$$
$$\begin{array}{c}	
	\mathscr P^5_{1313}=\frac{36b^2m^2r^5B}{B_1^4}=\frac{1}{\sin^2\theta}\mathscr P^5_{1414}, \ \ 	
	
	\mathscr P^5_{2323}=-\frac{36b^2m^2r^5}{B_1^2B}=\frac{1}{\sin^2\theta}\mathscr P^5_{2424}.
\end{array}.$$

From the above components we get the following:
\begin{pr}\label{pr5}
	The HBH spacetime fulfills the curvature conditions \\ $(i)$ $R\cdot P=-\frac{mB_2}{B_1^2}Q(g,P)$   $(ii)$ 
	$P\cdot S=-\frac{mB_2}{B_1^2}Q(g,S)$ and $(iii)$ the general form of $P$-compatible tensor is given by	
	$$
	\left(
	\begin{array}{cccc}
		\mathscr{Z}_{11} &\mathscr{Z}_{12} & 0 & 0 \\
		\mathscr{Z}_{12} & \mathscr{Z}_{22} & 0 & 0 \\
		0 & 0 & \mathscr{Z}_{33} & \mathscr{Z}_{34} \\
		0 & 0 & \mathscr{Z}_{34} & \mathscr{Z}_{44}
	\end{array}
	\right)
	$$
	where $\mathscr{Z}_{ij}$ are arbitrary scalars.

\end{pr}
\indent The non-vanishing components   of the concircular curvature tensor $W$  and conharmonic curvature tensor $K$  (upto symmetry) are computed as follows:
$$\begin{array}{c}
	W_{1212}=\frac{2mr^3(13b^2m-r^3)}{B_1^3},\ \ 
	W_{1313}=-\frac{mr^5B_2B}{B_1^4}=\frac{1}{\sin^2\theta}W_{1414}, \\
	W_{2323}=\frac{mr^5B_2}{B_1^2 B}=\frac{1}{\sin^2\theta}W_{2424}, \ \ 
	W_{3434}=\frac{2mr^7(5b^2m+r^3)\sin^2\theta}{B_1^3};
\end{array}.$$

$$\begin{array}{c}
	K_{1212}=\frac{2mB_2}{B_1^2},\ \ 
	K_{1313}=-\frac{mr^2B_2^2B}{B_1^4}=\frac{1}{\sin^2\theta}K_{1414}, \\
	
	K_{2323}=\frac{mr^2B_2^2}{B_1^2 B}=\frac{1}{\sin^2\theta}K_{2424}, \ \ 
	K_{3434}=-\frac{2mr^4B_2\sin^2\theta}{B_1^2}.
\end{array}.$$

If $\mathcal{M}^5=W \cdot R$ and $\mathcal{M}^6=K \cdot R$, then the non-vanishing components  of the tensors $\mathcal{M}^5$ and $\mathcal{M}^6$ are given as below:
$$\begin{array}{c}
	\mathcal{M}^5_{1223,13}=-\frac{3m^2r^8(40b^4m^2-14b^2mr^3+r^6)}{B_1^6}=-\mathcal{M}^5_{1213,23}, \\
	\mathcal{M}^5_{1434,13}=-\frac{3m^2r^{10}BB_2\sin^2\theta}{B_1^6}=-\mathcal{M}^5_{1334,14}, \\
	
	\mathcal{M}^5_{1224,14}=-\frac{3m^2r^8(40b^4m^2-14b^2mr^3+r^6)\sin^2\theta}{B_1^6}=-\mathcal{M}^5_{1214,24}, \\
	\mathcal{M}^5_{2434,23}=-\frac{3m^2r^{10}B_2\sin^2\theta}{B_1^4B}=-\mathcal{M}^5_{2334,24}
\end{array}$$
	
$$\begin{array}{c}
	\mathcal{M}^6_{1223,13}=\frac{3m^2r^5B_4B_2^2}{B_1^6}=-\mathcal{M}^6_{1213,23}, \ \ 
	\mathcal{M}^6_{1434,13}=-\frac{3m^2r^7B_2^2B \sin^2\theta}{B_1^6}=-\mathcal{M}^6_{1334,14}, \\
	\mathcal{M}^6_{1224,14}=\frac{3m^2r^5B_4B_2^2\sin^2\theta}{B_1^6}=-\mathcal{M}^6_{1214,24}, \ \ 
	
	\mathcal{M}^6_{2434,23}=\frac{3m^2r^7B_2^2\sin^2\theta}{B_1^4 B}=-\mathcal{M}^6_{2334,24}. 
\end{array}$$

From the above components we get the following:
\begin{pr}\label{pr6}
	The HBH spacetime fulfills the curvature conditions
	$$ (i)\ W\cdot R=-\frac{mr^3B_2}{B_1^3}Q(g,R),$$ \
	$(ii)$ 
	$K\cdot R=\frac{mB_2^2}{B_1^3}Q(g,R)$ and $(iii)$ the general form of $W$and $K$-compatible tensor are given by 	$$
	\left(
	\begin{array}{cccc}
	\mathscr{Z}_{11} &\mathscr{Z}_{12} & 0 & 0 \\
	\mathscr{Z}_{12} & \mathscr{Z}_{22} & 0 & 0 \\
	0 & 0 & \mathscr{Z}_{33} & \mathscr{Z}_{34} \\
	0 & 0 & \mathscr{Z}_{34} & \mathscr{Z}_{44}
	\end{array}
	\right)
	$$
	where $\mathscr{Z}_{ij}$ are arbitrary scalars.

\end{pr}

From  the above propositions $(\ref{pr1})$-$(\ref{pr6})$,  we can conclude  the curvature restricted  geometric properties of  HBH spacetime as follows:
\begin{thm}
	The HBH  spacetime admits the following curvature properties:
	\begin{enumerate}[label=(\roman*)]
		
		\item  $R\cdot R=-\frac{mB_2}{B_1^2} Q(g,R).$  Hence $R\cdot S=-\frac{mB_2}{B_1^2}Q(g,S)$, $R\cdot C=-\frac{mB_2}{B_1^2}Q(g,C)$, $R\cdot P=-\frac{mB_2}{B_1^2}Q(g,P)$, $R\cdot W=-\frac{mB_2}{B_1^2}Q(g,W)$ and $R\cdot K=-\frac{mB_2}{B_1^2}Q(g,K)$;
		%----------------------------------------------------------
		\item $C\cdot R=-\frac{mr^3B_2}{B_1^3}Q(g,R).$  Hence $C\cdot S=-\frac{mr^3B_2}{B_1^3}Q(g,S)$, $C\cdot C=-\frac{mr^3B_2}{B_1^3}Q(g,C)$, $C\cdot P=-\frac{mr^3B_2}{B_1^3}Q(g,P)$, $C\cdot W=-\frac{mr^3B_2}{B_1^3}Q(g,W)$ and $C\cdot K=-\frac{mr^3B_2}{B_1^3}Q(g,K)$;  
		%----------------------------------------------------------
		\item $W\cdot R=-\frac{mr^3B_2}{B_1^3}Q(g,R).$  Hence $W\cdot S=-\frac{mr^3B_2}{B_1^3}Q(g,S)$,\\ $W\cdot C=-\frac{mr^3B_2}{B_1^3}Q(g,C)$, $W\cdot P=-\frac{mr^3B_2}{B_1^3}Q(g,P)$, $W\cdot W=-\frac{mr^3B_2}{B_1^3}Q(g,W)$ and $W\cdot K=-\frac{mr^3B_2}{B_1^3}Q(g,K)$;
		 
		\item $K\cdot R=\frac{mB_2^2}{B_1^3}Q(g,R)$ . Hence $K\cdot S=\frac{mB_2^2}{B_1^3}Q(g,S)$, $K\cdot C=\frac{mB_2^2}{B_1^3}Q(g,C)$, $K\cdot P=\frac{mB_2^2}{B_1^3}Q(g,P)$, $K\cdot W=\frac{mB_2^2}{B_1^3}Q(g,W)$ and $K\cdot K=\frac{mB_2^2}{B_1^3}Q(g,K)$;
		  
		\item it  satisfies the pseudosymmetric type curvature conditions  $R\cdot R - \mathcal{L} Q(g,C) = Q(S,R) $, where $\mathcal{L}= \frac{m(16b^2m-r^3)}{B_2B_1}$;
		%------------------------------------------------------------
		
		\item the tensor $C\cdot R-R\cdot C$  depends linearly on  the tensors $Q(g,C)$, $Q(S,C)$, $Q(g,R)$ and $Q(S,R)$,
		 
		\item it is Ricci pseudosymmetric due to projective curvature i.e.,  $P\cdot S=-\frac{mB_2}{B_1^2}Q(g,S)$ is satisfied,
		%-----------------------------------------------------------
		
		%-----------------------------------------------------------
		\item  its conformal $2$-forms are recurrent for the for $1$-forms  $\{ 0, -\frac{6b^2m(8b^2m-5r^3)}{8b^4m^2r+2b^2mr^4-r^7},0 ,0 \},$ 
		%----------------------------------------------------------
		\item  it is a Roter type spacetime,
		%---------------------------------------------------------
		\item it is Ein$(2)$ spacetime as it possesses $S^2+\psi_1 S+\psi_2 g=0 $ for
		$\psi_1=\frac{12b^2m^2B_2}{B_1^3}$ and $\psi_2=\frac{288b^4m^4B_3}{B_1^5},$
		\item it is $2$-quasi-Einstein as $\alpha=-\frac{12b^2m^2}{B_1^2 },$
		
		\item  for $\alpha=-\frac{12b^2m^2}{B_1^2}$,  $\beta=1$, $\gamma=1$,    $\phi$= $\left\lbrace \frac{36b^2m^2r^3+B_1^2B}{2B_1^3},\frac{18b^2m^2r^3}{B_1^2B }-\frac{1}{2},0,0\right\rbrace $    and $\Pi$=$\left\lbrace -\frac{B}{B_1},1,0,0\right\rbrace $,  the HBH spacetime is generalized quasi-Einstein in the sense of Chaki,
		 
		\item the general form of $R$, $C$, $P$, $W$ and $K$-compatible tensors in HBH spacetime are given by 
		
		$$
		\left(
		\begin{array}{cccc}
			\mathscr{Z}_{11} &\mathscr{Z}_{12} & 0 & 0 \\
			\mathscr{Z}_{12} & \mathscr{Z}_{22} & 0 & 0 \\
			0 & 0 & \mathscr{Z}_{33} & \mathscr{Z}_{34} \\
			0 & 0 & \mathscr{Z}_{34} & \mathscr{Z}_{44}
		\end{array}
		\right)
		$$
		where $\mathscr{Z}_{ij}$ are arbitrary scalars,

		%\begin{equation}
		%	g^{\mu\nu}=
		%	\begin{pmatrix}
		%	$$	\mathcal{T}_{11} & \mathcal{T}_{12} & 0 & 0  \\
		%		\mathcal{T}_{21} & \mathcal{T}_{22} &0 & 0\\
		%		0 & 0 & \mathcal{T}_{33} & \mathcal{T}_{34}\\
		%		0 & 0 &\mathcal{T}_{34} & \mathcal{T}_{44}$$
		%	\end{pmatrix}\quad \quad
		%g^{\mu\nu}=\begin{pmatrix}
		%-1 & 0 & 0 & 0  \\
		%0 & \alpha^2 &0 & 0\\
		%0 & 0 & \frac{1}{r^2} & 0\\
		%0 & 0 & 0 & \frac{1}{r^2 \sin^2 \theta}.
		%\end{pmatrix}
		%\label{2}
		%\end{equation}

		\item its Ricci tensor is compatible for $C$, $P$, $R$, $K$ and $W$.
	\end{enumerate}
	
\end{thm}
%The non-zero components of  the tensor $R\cdot S$ are given by 
%$$\begin{array}{c}
%R\cdot S_{1313}=-\frac{36b^2m^3r^5(4b^2m-r^3)(2b^2m+(-2m+r)r^2)}{(2b^2m+r^3)^6}; \\

%R\cdot S_{1414}=-\frac{36b^2m^3r^5(4b^2m-r^3)(2b^2m+(-2m+r)r^2)\sin^2\theta}{(2b^2m+r^3)^6}; \\

%R\cdot S_{2323}=\frac{36b^2m^3r^5(4b^2m-r^3)}{(2b^2m+r^3)^4(2b^2m+(-2m+r)r^2)}; \\
%R\cdot S_{2424}=\frac{36b^2m^3r^5(4b^2m-r^3)\sin^2\theta}{(2b^2m+r^3)^4(2b^2m+(-2m+r)r^2)}; \\

%\end{array}$$

%The non-zero components of the Tachibana tensor $Q(g,S)$ are given by
%$$\begin{array}{c}
%Q(g,S)_{1313}=\frac{36b^2m^2r^5(2b^2m+r^2(-2m+r))}{(2b^2m+r^3)^4}; \\
%Q(g,S)_{1414}=\frac{36b^2m^2r^5(2b^2m+r^2(-2m+r))\sin^2\theta}{(2b^2m+r^3)^4}; \\
%Q(g,S)_{2323}=-\frac{36b^2m^2r^5}{(2b^2m+r^3)^2(2b^2m+r^2(-2m+r))}; \\
%Q(g,S)_{2424}=-\frac{36b^2m^2r^5\sin^2\theta}{(2b^2m+r^3)^2(2b^2m+r^2(-2m+r))}; \\

%\end{array}$$

\begin{rem}
	The HBH spacetime does not admit the following geometric structures:
	\begin{enumerate}[label=(\roman*)]
		\item $\nabla P\neq 0$ and hence $\nabla R\neq 0$, $\nabla C\neq 0$, $\nabla K\neq 0$ and $\nabla W\neq 0$,
		
		\item  for any $1$-form $\Pi$, $\nabla P\neq \Pi \otimes P$ and hence it is not recurrent for $P$, $R$, $W$, $K$ and $C$,
		
		\item it does not satisfy the semi-symmetric type condition $R\cdot H =0$ where $H=P, K, W, C, S$,
		
		\item it is not Ricci generalized pseudosymmetric, 
		
		\item it does not realize $P\cdot R= \mathcal{L}Q(g,R)$ for any smooth function $\mathcal{L}$. Hence it is neither   $P\cdot W= \mathcal{L}Q(g,W)$, $P\cdot K= \mathcal{L}Q(g,K)$ nor $P\cdot C= \mathcal{L}Q(g,C),$
		
		\item it is not $T$-space by Venzi  for $T=C, R, P, W, K,$ 
		
		\item it is neither Einstein nor quasi-Einstein,
		\item the curvature $2$-forms  for $R$, $K$, $W$ and $P$ are not recurrent, 
		
		\item the Ricci tensor of HBH spacetime is neither cyclic parallel nor Codazzi type,
		\item the HBH spacetime is neither  weakly symmetric nor Chaki pseudosymmetric for $P,$ $W,$ $K,$ $R$ and $C$. 
		
	\end{enumerate}  
\end{rem}

\section{\bf Energy momentum tensor of Hayward black hole spacetime}
%%%%%%%%%%%%%%%%%%%%%%%%%%%%%%%%%%%%%%%%%%%%%%%%%%%%%%%%%%%%%%%%%%%%%%%
%=======================================================================

In Einstein field equation (briefly, EFE), the energy  momentum tensor $T^{EM}$ in terms of curvature restrictions  is presented as 
$T^{EM}=\frac{1}{\nu}[S-\frac{\kappa}{2}g+\Lambda g],$
where $\Lambda$ is the cosmological constant, $\nu=\frac{8\pi G}{c^4}$ ($G$ being the Newton's gravitational constant and c being the speed of light in vacuum).
The components other than zero of the energy momentum tensor $T^{EM}_{ab}$ are  given below: 

$$\begin{array}{c}
	T^{EM}_{11}=-\frac{3b^2m^2B}{2B_1^3}, \ \
	T^{EM}_{22}=-\frac{3b^2m^2}{2B_1B}, \\
	T^{EM}_{33}=\frac{3b^2m^2r^2B_3}{B_1^3}=\frac{1}{\sin^2\theta}T^{EM}_{44}. \\

\end{array}$$
The non-vanishing components of the covariant derivative of  energy momentum tensor are calculated as  follows:
$$\begin{array}{c}
	T^{EM}_{11,2}=\frac{9b^2m^2r^2B}{B_1^4}, \ \ 
	T^{EM}_{22,2}=-\frac{9b^2m^2r^2}{B_1^2B}, \\
	T^{EM}_{23,3}=\frac{9b^2m^2r^4}{2B_1^3}=\frac{1}{\sin^2\theta}T^{EM}_{24,4}, \\

	T^{EM}_{33,2}=\frac{9b^2m^2r^4(-5b^2m+2r^3)}{B_1^4}=\frac{1}{\sin^2\theta}T^{EM}_{44,2}.	
\end{array}$$
Let $\mathcal{V}^1=R \cdot T^{EM}$, $\mathcal{V}^2=C \cdot T^{EM}$, $\mathcal{V}^3=W \cdot T^{EM}$, $\mathcal{V}^4=K \cdot T^{EM}$, and $\mathcal{U}^1=Q(g,T^{EM})$.  Then the non-vanishing components   of the tensors $\mathcal{V}^1$, $\mathcal{V}^2$, $\mathcal{V}^3$, $\mathcal{V}^4$ and $\mathcal{U}^1$ are obtained as follows:
$$\begin{array}{c}
	\mathcal{V}^1_{1313}=-\frac{9b^2m^3r^5B_2B}{2B_1^6}=\frac{1}{\sin^2\theta}\mathcal{V}^1_{1414}, \ \ 	
	\mathcal{V}^1_{2323}=\frac{9b^2m^3r^5B_2}{2B_1^4B}=\frac{1}{\sin^2\theta}\mathcal{V}^1_{2424}; 	
\end{array}$$
$$\begin{array}{c}
	\mathcal{V}^2_{1313}=-\frac{9b^2m^3r^8B_2B}{2B_1^7}=\frac{1}{\sin^2\theta}\mathcal{V}^2_{1414}, \ \	
	\mathcal{V}^2_{2323}=\frac{9b^2m^3r^8B_2}{2B_1^5B}=\frac{1}{\sin^2\theta}\mathcal{V}^2_{2424}; 
\end{array}$$
$$\begin{array}{c}
	\mathcal{V}^3_{1313}=-\frac{9b^2m^3r^8B_2B}{2B_1^7}=\frac{1}{\sin^2\theta}\mathcal{V}^3_{1414}, \ \	
	\mathcal{V}^3_{2323}=\frac{9b^2m^3r^8B_2}{2B_1^5B}=\frac{1}{\sin^2\theta}\mathcal{V}^3_{2424}; 
\end{array}$$
$$\begin{array}{c}
	\mathcal{V}^4_{1313}=-\frac{9b^2m^3r^5B_2^2B}{2B_1^7}=\frac{1}{\sin^2\theta}\mathcal{V}^4_{1414}, \ \	
	\mathcal{V}^4_{2323}=-\frac{9b^2m^3r^5B_2^2}{2B_1^5B}=\frac{1}{\sin^2\theta}\mathcal{V}^4_{2424}; 
\end{array}$$
$$\begin{array}{c}
	\mathcal{U}^1_{1313}=\frac{9b^2m^2r^5B}{2B_1^4}=\frac{1}{\sin^2\theta}\mathcal{U}^1_{1414}, \ \ 
	\mathcal{U}^1_{2323}=-\frac{9b^2m^2r^5}{2B_1^2B}=\frac{1}{\sin^2\theta}\mathcal{U}^1_{2424}.
\end{array}$$
\indent From the above components we get the following theorem:
\begin{thm}
	The energy momentum tensor of the HBH spacetime admits the following geometric properties:  
	\begin{enumerate}[label=(\roman*)]
		\item  $R\cdot  T^{EM}=-\frac{mB_2}{B_1^2}Q(g,T^{EM})$ i.e., the nature of the energy momentum tensor is pseudosymmetric, 
		\item  $C\cdot  T^{EM}=-\frac{mr^3B_2}{B_1^3}Q(g,T^{EM})$,
		
		\item $W\cdot T^{EM}=-\frac{mr^3B_2}{B_1^3}Q(g,T^{EM}),$ 
		\item $K\cdot T^{EM}=\frac{mB_2^2}{B_1^3}Q(g,T^{EM})$  and

		\item the energy momentum tensor is Riemann compatible, projective compatible, conharmonic compatible, concircular compatible and conformal compatible.	 
	\end{enumerate}
\end{thm}

%=============================================================================
%
%%%%%%%%%%%%%%%%%%%%%%%%%%%%%%%%%%%%%%%%%%%%%%%%%%%%%%%%%%%%%%%%%%%%%%%%%%

%%%%%%%%%%%%%%%%%%%%%%%%%%%%%%%%%%%%%%%%%%%%%%%%%%%%%%%%%%%%%%%%%%%%%%%
%=======================================================================
% 
\section{\bf Hayward black hole spacetime Vs interior black hole spacetime and Reissner-Nordstr\"om spacetime  }
%%%%%%%%%%%%%%%%%%%%%%%%%%%%%%%%%%%%%%%%%%%%%%%%%%%%%%%%%%%%%%%%%%%%%%%
%=======================================================================
The interior black hole spacetime \cite{Doran_interior_2008,SDHK_interior_2020}  is a spherically symmetric non-static solution of EFE. Physically, it describes the empty spacetime in the exterior region of a black hole. A comparison between HBH spacetime and interior black hole spacetime in terms of their curvature properties is delineated as follows:\\
\noindent\textbf{Similarities:}
\begin{enumerate}[label=(\roman*)]
	\item both the spacetimes are pseudosymmetric,
	%-----------------------------------------------------------------
	
	\item both the spacetimes are pseudosymmetric due to conharmonic, concircular as well as conformal  curvature, 
	%---------------------------------------------------------
	\item both the spacetimes are  Einstein manifolds of level $2$ and $2$-quasi Einstein manifolds,
	%------------------------------------------------------------------
	\item both the spacetimes are Roter type,
	%----------------------------------------------------------
	
	\item Ricci tensor is  Riemann compatible as well as Weyl compatible .
	
\end{enumerate}

\indent Again, the exterior gravitational field of a non-rotating charged body can be described by Reissner-Nordstr\"om spacetime \cite{Kowa06}, which  is a spherically symmetric solution of EFE having  cosmological constant zero. This solution is more general than the Schwarzschild solution of EFE as the Reissner-Nordstr\"om solution admits non-vanishing charges. An elegant comparison between HBH spacetime and Reissner-Nordstr\"om spacetime based on the curvature properties is described as follows: \\  
\noindent\textbf{Dissimilarities:}
\begin{enumerate}[label=(\roman*)]
	\item the conharmonic $2$-forms of Reissner-Nordstr\"om spacetime  are recurrent while HBH spacetime does not admit such recurrence,
	
	\item HBH spacetime does not vanish scalar curvature while for the  Reissner-Nordstr\"om spacetime  the scalar curvature vanishes.
	%----------------------------------------------------------
	%\item the Bardeen spacetime comes out with a weakly generalized recurrent manifold while  Reissner-Nordstr\"om spacetime doesn't,
	%------------------------------------------------------------
	%\item also the Bardeen spacetime admits special recurrent like structure $\nabla R=A\otimes (g\wedge S)$ ($A$ being some $1$-form) but  Reissner-Nordstr\"om spacetime does not admit such recurrence.
	
\end{enumerate}

However, the HBH spacetime and the Reissner-Nordstr\"om spacetime have the following similar properties:
%\noindent\textbf{Similarities:}
\begin{enumerate}[label=(\roman*)]
	\item both spacetimes are Roter type,
	%-----------------------------------------------------------------
	\item both the spacetimes are  Einstein manifolds of level $2$,
	%------------------------------------------------------------------
	\item both are pseudosymmetric as well as pseudosymmetric due to Weyl conformal tensor,
	
	%---------------------------------------------------------
	\item conformal $2$-forms for both the spacetimes are  recurrent,
	%----------------------------------------------------------
	\item both are $2$-quasi-Einstein manifold,
	\item  Ricci tensor of both the spacetimes are  Riemann compatible as well as Weyl compatible.
	
\end{enumerate} 

\section{\bf Ricci soliton and symmetries on Hayward black hole spacetime}

Let $\mathcal{K}(M)$ be the set of all Killing vector fields on $M$. Then $\mathcal{K}(M)$ is a Lie subalgebra of the Lie algebra $\chi(M)$ of all smooth vector fields on $M$ and $\mathcal{K}(M)$ contains at most $n(n+1)\slash 2$ linearly independent Killing vector fields, and if $\mathcal{K}(M)$ consists of exactly $n(n+1)\slash 2$ linearly independent Killing vector fields, then $M$ is known as a maximally symmetric space. We mention that $M$ is a maximally symmetric space if $M$ is of constant scalar curvature. We note that the scalar curvature $\kappa$ of HBH spacetime is not constant as shown in Section 3 by $\kappa=\frac{24b^2m^2(r^3-4b^2m)}{(2b^2m+r^3)^3}$ and hence it is not maximally symmetric. Now, we investigate some Killing and non-Killing vector fields on HBH spacetime given as follows:

\begin{pr}
	The vector fields ${\frac{\partial}{\partial t}}$ and ${\frac{\partial}{\partial \phi}}$ on the HBH spacetime are Killing    (i.e., $\pounds_{\frac{\partial}{\partial t}}g=0=\pounds_{\frac{\partial}{\partial \phi}}g$).
\end{pr}

\begin{cor}
	For each real number $\lambda_1$ and $\lambda_2$, the vector field $\lambda_1{\frac{\partial}{\partial t}}+\lambda_2{\frac{\partial}{\partial \phi}}$ on the HBH spacetime is also Killing.
\end{cor}

The vector field ${\frac{\partial}{\partial r}}$ is non-Killing, and if $\mathcal A=\pounds_{\frac{\partial}{\partial r}}g$, then the non-zero components of   $\mathcal A$ are calculated as follows:
$$\begin{array}{c}
	\mathcal A_{11}=\frac{2mrB_2}{B_1^2},\ \ \ 
	\mathcal A_{22}=\frac{2mrB_2}{B^2},%{\{2b^2m+r^2(r-2m)\}^2},
	\nonumber\\
	\mathcal A_{33}=2r,\ \ \ 
	\mathcal A_{44}=2r\sin^2\theta.\nonumber
\end{array}$$
Therefore, for the non-Killing vector field $\frac{\partial}{\partial r}$ and the $1$-form $\eta=(0,1,0,0)$, the HBH spacetime possesses the following relation:
\begin{eqnarray}\nonumber
	\pounds_{\frac{\partial}{\partial r}}g+ 2\sigma_1 S+2\sigma_2 g-2\sigma_3 \eta\otimes\eta=0,
\end{eqnarray}
where $\sigma_1,\sigma_2,\sigma_3$ are given by
\begin{equation}\label{LieCoefficient}
	\left.
	\begin{aligned}
		\sigma_1&=\frac{B_1^2(4b^4m^2+4b^2mr^3-3mr^5+r^6)}{36b^2m^2r^4B},\\
		\sigma_2&=\frac{4b^4m^2-2b^2mr^3+3mr^5-2r^6}{3r^4B},\\
		\sigma_3&=\frac{2mrB_2}{B^2}.
	\end{aligned}\ \ 
	\right\rbrace	
\end{equation}

This leads to the following:

\begin{thm}
	The HBH spacetime  realizes  almost $\eta$-Ricci-Yamabe soliton for the non-Killing soliton vector field $\frac{\partial}{\partial r}$ and the 1-form $\eta=(0,1,0,0)$ provided $(2b^2m-2mr^2+r^3)\neq0$, i.e., for the soliton vector field $\xi=\frac{\partial}{\partial r}$, the HBH spacetime  possesses
	$$\frac{1}{2}\pounds_\xi g+\sigma_1 S+\left(\lambda-\frac{1}{2}\sigma_4\kappa\right) g+\sigma_3 \eta\otimes\eta=0,$$
	where $\sigma_4=2$,   $\lambda=\sigma_2+\kappa$, and $\sigma_1,\sigma_2,\sigma_3$ are given in (\ref{LieCoefficient}).

\end{thm}

\begin{thm}If $(2b^2m+r^3)^2(4b^4m^2+4b^2mr^3-3mr^5+r^6)=36b^2m^2r^4(2b^2m-2mr^2+r^3)$ with $(2b^2m-2mr^2+r^3)\neq0$, then
	for the soliton vector field $\frac{\partial}{\partial r}$, the HBH spacetime  admits an almost $\eta$-Ricci soliton  with the $1$-form $\eta=(0,1,0,0)$, i.e., for the vector field $\xi=\frac{\partial}{\partial r}$, the HBH spacetime realizes
	\begin{eqnarray}
		&&\frac{1}{2}\pounds_{\xi}g+S+\sigma_2 g-\sigma_3\eta\otimes\eta=0.\nonumber
	\end{eqnarray}
where $\sigma_2$, $\sigma_3$ is given in (\ref{LieCoefficient}).
\end{thm}

Let $\mathcal E=\pounds_{\frac{\partial}{\partial r}}S$, $\mathcal G=\pounds_{\frac{\partial}{\partial r}}\widetilde{R}$ and $\mathcal H=\pounds_{\frac{\partial}{\partial r}}R$. Then the non-vanishing components of $\mathcal E$, $\mathcal G$ and $\mathcal H$ are computed as follows:
$$\begin{array}{c}
 \mathcal E_{11}=-\frac{24b^2m^2r^2\{2(7m-3r)r^6+b^2mr^3(3r-40m)+2b^4m^2(4m+15r) \} }{B_1^5},\\
	\mathcal E_{22}=-\frac{24b^2m^2r\{2b^4m^2(4m-15r)+b^2m(20m-3r)r^3+2r^6(3r-5m)  \} }{B_1^3 B^2 },\\
	\mathcal E_{33}=-\frac{48b^2m^2rB_3}{B_1^3}=\frac{1}{\sin^2\theta} \mathcal E_{44}
\end{array}$$

$$\begin{array}{c}
	\mathcal G^1_{212}=\frac{2mr\left\lbrace16b^6m^3(2m-15r)+r^9(4m-3r)+24 b^4m^2r^3(5m+3r)+6b^2mr^6(15r-26m)\right\rbrace}{B_1^3B} \\ %\left\lbrace 2b^2m+r^2(r-2m)\right\rbrace^2} \\
	\mathcal G^1_{313}=\frac{mr(16b^4m^2-26b^2mr^3+r^6)}{B_1^3}=-\mathcal G^1_{331}=-\frac{1}{\sin^2\theta}\mathcal G^1_{441}=\mathcal G^2_{323} \\=-\mathcal G^1_{332}=\frac{1}{\sin^2\theta}\mathcal G^2_{424}=-\frac{1}{\sin^2\theta}\mathcal G^2_{442}, \\
	\mathcal G^2_{112}=\frac{2mr\left\lbrace(8m-3r)r^9+72b^4m^2r^3(5m+r)+6b^2mr^6(15r-38m)-16b^6m^3(2m+15r)\right\rbrace}{B_1^5}=-\mathcal G^2_{121}, \\
	
	\mathcal G^3_{113}=-\frac{mr\left\lbrace(8m-3r)r^6+4b^2mr^3(6r-19m)+4b^4m^2(8m+15r)\right\rbrace}{B_1^4}=-\mathcal G^3_{131}=\mathcal G^4_{114}=-\mathcal G^4_{141}, \\
	
	\mathcal G^3_{223}=\frac{mr\left\lbrace(4m-3r)r^6+4b^2mr^3(6r-11m)+4b^4m^2(15r-8m)\right\rbrace}{B_1^2 B^2}=-\mathcal G^3_{232}=\mathcal G^4_{224}=-\mathcal G^4_{242}, \\
	\mathcal G^3_{434}=\frac{2mrB_2\sin^2\theta}{B_1^2}=-\mathcal G^3_{443}=-\sin^2\theta \mathcal G^4_{334}=\sin^2\theta \mathcal G^4_{343}.
\end{array}$$

$$\begin{array}{c}
	\mathcal H_{1212}=\frac{6m(40b^4m^2r^2-32b^2mr^5+r^8)}{B_1^4}=-\mathcal H_{1221}=\mathcal H_{2121}, \\
	\mathcal H_{1313}=\frac{mr\left\lbrace -32b^6m^3+(4m-r)r^8+4b^2mr^5(6r-17m)+4b^4m^2r^2(16m+9r) \right\rbrace}{B_1^4}=\mathcal H_{3131}=-\mathcal H_{1331} \\ =\frac{1}{\sin^2\theta} \mathcal H_{4141}=\frac{1}{\sin^2\theta} \mathcal H_{1414}=-\frac{1}{\sin^2\theta} \mathcal H_{1441}, \\
	\mathcal H_{2323}=\frac{mr\left\lbrace 32b^6m^3-36b^4m^2r^3+12b^2m(3m-2r)r^5+r^9 \right\rbrace}{B_1^2 B^2 }=-\mathcal H_{2332}=\frac{1}{\sin^2\theta} \mathcal H_{2424} \\=-\frac{1}{\sin^2\theta} \mathcal H_{2442}= \mathcal H_{3232}=\frac{1}{\sin^2\theta} \mathcal H_{4242}, \\
	\mathcal H_{3434}=\frac{2mr^3(8b^2m+r^3)\sin^2\theta}{B_1^2}=-\mathcal H_{3443}=\mathcal H_{4343}.
\end{array}$$

If $\mathcal M=\pounds_{\frac{\partial}{\partial \theta}}g$, $\mathcal N=\pounds_{\frac{\partial}{\partial \theta}}S$, $\mathcal Q=\pounds_{\frac{\partial}{\partial \theta}}\widetilde{R}$ and $\mathcal O=\pounds_{\frac{\partial}{\partial \theta}}R$, then the non-zero components of $\mathcal M$, $\mathcal N$, $\mathcal Q$ and $\mathcal O$ are given as follows:
$$\begin{array}{c}
	\mathcal Q^1_{414}=\frac{mr^2B_2\sin2\theta}{B_1^2}=-\mathcal Q^1_{441}=\mathcal Q^2_{424}=-\mathcal Q^2_{442}, \ \
	\mathcal Q^3_{434}=\frac{2mr^2\sin2\theta}{B_1}=-\mathcal Q^3_{443},
\end{array}$$

$$\begin{array}{c}
		\mathcal O_{1414}=\frac{mr^2B \sin2\theta}{B_1^3}=-\mathcal O_{1441}=\frac{1}{B_2}\mathcal O_{4141}, \ \
	\mathcal O_{2424}=\frac{mr^2B_2\sin2\theta}{B_1B}=-\mathcal O_{2442}=\mathcal O_{4242}, \\
	\mathcal O_{3434}=\frac{2mr^4\sin2\theta}{B_1}=-\mathcal O_{3443}=\mathcal O_{1343}.
\end{array}$$

From the above calculation of the Lie derivative of various curvature tensors it can be easily checked that with respect to the non-Killing vector fields $\frac{\partial}{\partial r}$, $\frac{\partial}{\partial \theta}$ and  $\lambda_1\frac{\partial}{\partial r}+\lambda_2\frac{\partial}{\partial \theta}$ ($\lambda_1,\lambda_2$ being real numbers), the HBH spacetime  admits 
	 \begin{enumerate}[label=(\roman*)]
	 	\item neither Ricci collineation nor Ricci inheritance,
	 	
	 	\item neither curvature collineation for (1,3)-type curvature tensor nor curvature collineation for (0,4)-type curvature tensor,
	 	
	 	\item neither curvature inheritance for (1,3)-type curvature tensor nor curvature inheritance for (0,4)-type curvature tensor.
	 \end{enumerate}

\section{\bf Hayward black hole spacetime Vs point-like global monopole spacetime  }

The point-like global monopole spacetime \cite{Barriola1989,SAD_pgm_2023} is a static and spherically symmetric solution of EFE. It is a heavy object characterized by divergent mass and spherically symmetry, and against polar as well as spherical perturbation it is expected to be stable. A comparative study between the HBH spacetime and the point-like global monopole spacetime with respect to various kind of symmetries and Ricci soliton is given as follows:\\
\noindent\textbf{Similarities:}
\begin{enumerate}[label=(\roman*)]
	\item both the spacetimes admits motion for the vector fields $\frac{\partial}{\partial t}$ and $\frac{\partial}{\partial \phi}$, i.e.,  the vector fields $\frac{\partial}{\partial t}$ and $\frac{\partial}{\partial \phi}$ are Killing in both the spacetimes,
	
	\item the vector fields $\frac{\partial}{\partial r}$ and $\frac{\partial}{\partial \theta}$ are non-Killing in both the spacetimes,
	
	\item with respect to the non-Killing vector field $\frac{\partial}{\partial \theta}$, both the spacetimes realize neither curvature collineation nor curvature inheritance for (1,3)-type curvature tensor,
	
	\item  with respect to the non-Killing vector field $\frac{\partial}{\partial \theta}$, both the spacetimes possess neither Ricci collineation nor Ricci inheritance.
\end{enumerate}
Nevertheless, they have the following dissimilar properties:\\
\noindent\textbf{Dissimilarities:}

\begin{enumerate}[label=(\roman*)]
	\item with respect to the non-Killing vector field $\frac{\partial}{\partial r}$, the point-like global monopole spacetime admits Ricci collineation as well as curvature collineation for (1,3)-type curvature tensor, whereas HBH spacetime does not admit such collineations,
	
	\item for the non-Killing vector fields $\frac{\partial}{\partial r}$, $\frac{\partial}{\partial \theta}$ and  $\lambda_1\frac{\partial}{\partial r}+\lambda_2\frac{\partial}{\partial \theta}$ ($\lambda_1,\lambda_2$ being real numbers), the point-like global monopole spacetime possesses curvature inheritance for the (0,4)-type curvature tensor, but HBH spacetime does not realize such inheritance,
	
	\item with respect to the soliton vector field $\frac{\partial}{\partial r}$, the HBH spacetime admits both the almost $\eta$-Ricci soliton and almost $\eta$-Ricci-Yamabe soliton for the 1-form $\eta=(0,1,0,0)$, but the point-like global monopole spacetime  realizes neither almost $\eta$-Ricci soliton nor almost $\eta$-Ricci-Yamabe soliton with respect to the non-Killing vector field $\frac{\partial}{\partial r}$.
\end{enumerate}

%%%%%%%%%%%%%%%%%%%%%%%%%%%%%%%%%%%%%%%%%%%%%%%%%%%%%%%%%%%%%%%%%%%%%%%%%%%%%%%%%%%%%%%%%%%%%%%%%%%%%%%%%%%%%%%%
\section{\bf  Acknowledgment}
The third author is grateful to the Council of Scientific and Industrial Research (CSIR File No.: 09/025(0253)/2018-EMR-I), Govt. of India, for the award of SRF (Senior Research Fellowship). The fourth author greatly acknowledges to The University Grants Commission, Government of India for the award of Senior Research Fellow. All the algebraic computations of Section $3$ to $7$ are performed by a program in Wolfram Mathematica developed by the first author (A. A. Shaikh).

%%%%%%%%%%%%%%%%%%%%%%%%%%%%%%%%%%%%%%%%%%%%%%%%%%%%%%%%%%%%%%%%%%%%%%%%%%%%%%%%%%%%%%%%%%%%%%%%%%%%%%%%%%%%%%%%

%
%%%%%%%%%%%%%%%%%%%%%%%%%%%%%%%%%%%%%%%%%%%%%%%%%%%%%%%%%%%%%%%%%%%%%%%%%%%%%%%%%%%%%%%%%%%%%%%%%%%%%%%%%%%%%%%%%
%%%%%%%%%%%%%%%%%%%%%%%%%%%%%%%%%%%%%%%%%%%%%%%%%%%%%%%%%%%%%%%%%%%%%%%%%%%%%%%%%%%%%%%%%%%%%%%%%%%%%%%%%%%%%%%%%

%%%%%%%%%%%%%%%%%%%%%

%%%%%%%%%%%%%%%%%%%%%%%%%%%%%%%%%%%%%%%%%%%%%%%%%%%%%%%%%%%%%%%%%%%%%%%%%%%%%%%%%%%%%%%%%%%%%%%%%%%%
\end{document}